\documentclass[12pt]{amsart}
\usepackage{amsmath}
\usepackage{amscd}
\usepackage[psamsfonts]{amssymb}
\usepackage[mathscr,psamsfonts]{eucal}
\DeclareMathAlphabet{\eurm}{U}{eur}{m}{n}
\DeclareMathAlphabet{\cyrm}{U}{UWCyr}{m}{n}
\newcommand{\M}[1]{{\mathscr{#1}}}

\newcommand{\E}[1]{{\eurm{#1}}}

\newcommand{\f}[1]{{\boldsymbol{#1}}}

\DeclareMathOperator{\id}{{id}}
\DeclareMathOperator{\Ob}{{Ob}}
\DeclareMathOperator{\Mor}{{Mor}}
\DeclareMathOperator{\proj}{{pr}}
\DeclareMathOperator{\lin}{{lin}}
\DeclareMathOperator{\Lin}{{Lin}}
\DeclareMathOperator{\Alt}{{Alt}}
\DeclareMathOperator{\pol}{pol}
\DeclareMathOperator{\inv}{inv}
\DeclareMathOperator{\Ker}{Ker}
\DeclareMathOperator{\Cla}{Cla}
\DeclareMathOperator{\pr}{pr}
\DeclareMathOperator{\byd}{\,{\raisebox{.1ex}{$\eurm :$}{\eurm =}}\,}
\newcommand{\sig}{\sigma}
\newcommand{\alp}{\alpha}
\newcommand{\bet}{\beta}
\newcommand{\lam}{\lambda}
\newcommand{\ome}{\omega}
\newcommand{\Lam}{\Lambda}

\newcommand{\gam}{\gamma}
\newcommand{\bEq}{\begin{eqnarray}}
\newcommand{\eEq}{\end{eqnarray}}
\newcommand{\beq}{\begin{eqnarray*}}
\newcommand{\eeq}{\end{eqnarray*}}

\newcommand{\ac}[1]{\acute{#1}}

\newcommand{\dt}[1]{{\dot{#1}}}
\newcommand{\uten}[1]{\underset{#1}{\otimes}}
\newcommand{\ucar}[1]{\underset{#1}{\times}}
\newcommand{\sep}[1]{\qquad\text{#1}\qquad}

\newcommand{\col}[3]{_{#1}{}^{#2}{}_{#3}}
\newcommand{\der}{\partial}
\newcommand{\nab}{\nabla}

\newcommand{\com}{\circ}
\newcommand{\car}{\times}
\newcommand{\ten}{\otimes}
\newcommand{\wed}{\wedge}
\newcommand{\Rn}{{I\!\!R}}

\newcommand{\sem}{\rtimes}

\begin{document}

\newcounter{theorem}

\newtheorem{definition}[theorem]{Definition}
\newtheorem{lemma}[theorem]{Lemma}
\newtheorem{proposition}[theorem]{Proposition}
\newtheorem{theorem}[theorem]{Theorem}
\newtheorem{corollary}[theorem]{Corollary}
\newtheorem{remark}[theorem]{Remark}
\newtheorem{example}[theorem]{Example}
\newtheorem{Note}[theorem]{Note}
\newcounter{assump}
\newtheorem{Assumption}{\indent Assumption}[assump]
\renewcommand{\thetheorem}{\thesection.\arabic{theorem}}

\newcommand{\bCr}{\begin{corollary}}
\newcommand{\eCr}{\end{corollary}}
\newcommand{\bDf}{\begin{definition}\em}
\newcommand{\eDf}{\end{definition}}
\newcommand{\bLm}{\begin{lemma}}
\newcommand{\eLm}{\end{lemma}}
\newcommand{\bPr}{\begin{proposition}}
\newcommand{\ePr}{\end{proposition}}
\newcommand{\bRm}{\begin{remark}\em}
\newcommand{\eRm}{\end{remark}}
\newcommand{\bEx}{\begin{example}\em}
\newcommand{\eEx}{\end{example}}
\newcommand{\bTh}{\begin{theorem}}
\newcommand{\eTh}{\end{theorem}}
\newcommand{\bNt}{\begin{Note}\em}
\newcommand{\eNt}{\end{Note}}
\newcommand{\bPf}{\begin{proof}[\noindent\indent{\sc Proof}]}
\newcommand{\ePf}{\end{proof}}

\title[Higher order reduction theorems]
        {Higher order reduction theorems
        \\for general linear connections}

\author{Josef Jany\v ska}



\address{
\newline
{Department of Mathematics, Masaryk University
\newline
Jan\'a\v ckovo n\'am. 2a, 662 95 Brno, Czech Republic}
\newline
E-mail: {\tt janyska@math.muni.cz}
}
\thanks{This paper has been supported
by the Ministry of Education of the Czech Republic under the Project
MSM 143100009.}

\begin{abstract}
The reduction theorems for general linear and classical connections
are generalized for operators with values
in higher order gauge-natural bundles.
We prove that natural operators depending on the $s_1$-jets of classical
connections, on the $s_2$-jets of general linear connections
and on the $r$-jets
of tensor fields with values in gauge-natural bundles of order
$k\ge 1$, $s_1+2\ge s_2$, $s_1,s_2\ge r-1\ge k-2$,
can be factorized through the $(k-2)$-jets of both
connections, the $(k-1)$-jets of the tensor fields and sufficiently high
covariant differentials of the curvature tensors and the tensor fields.
\end{abstract}
\maketitle

\section*{Introduction}
\setcounter{equation}{0}
\setcounter{theorem}{0}

It is well known that natural operators of classical  (linear
and symmetric)
connections on manifolds and of tensor fields with values
in natural bundles of order one
can be factorized through the curvature tensors, the tensor
fields and their covariant differentials. These theorems are known as
the first (operators on classical connections only)
and the second reduction, \cite{Lub72, Sch54}, or replacement,
\cite{Tho27, ThoMic27}, theorems.
In \cite{KolMicSlo93} the reduction theorems are proved by using
methods of natural bundles and operators, \cite{KruJan90, Nij72, Ter78}.

In \cite{Jan04} the reduction theorems were generalized for
general linear connections on vector bundles.
In this gauge-natural situation
we need auxiliary classical connections on the base manifolds.
It is proved that natural operators with values in gauge-natural
bundles of order (1,0)
defined
on the  space of general linear connections  on a
vector bundle, on the space of classical connections on
the base manifold and on certain tensor bundles can
be factorized through the curvature tensors
of linear and classical connections, the tensor
fields  and their
covariant differentials with respect to both connections.

In \cite{Jan04a} another generalization of the classical
reduction theorems was presented.
Namely, the reduction theorems were proved
 for operators with values in higher order
natural bundles. It was proved that an $r$-th order natural operator
on classical connections with values in natural bundles
of order $k\ge 1$, $r+2\ge k$,
can be factorized through the $(k-2)$-jets of connections and sufficiently
high covariant differentials of the curvature tensor.

In this paper we combine both possible generalizations of the
reduction theorems and we prove the reduction theorems for
general linear connections on vector bundles for operators with values
in higher order gauge-natural bundles.

All manifolds and maps are assumed to be smooth.
The sheaf of (local) sections of a fibered  manifold $p:\f Y\to \f X$
is denoted by $C^\infty(\f Y)$, $C^\infty(\f Y,\Rn)$ denotes the sheaf
of (local) functions.

\section{Gauge-natural bundles}
\setcounter{equation}{0}
\setcounter{theorem}{0}

Let $\M M_m$  be the category of $m$-dimensional
$C^\infty$-manifolds  and  smooth  embeddings. Let $\M F\M M_m$
be the category of smooth fibered  manifolds over  $m$-dimensional
bases and smooth fiber manifold maps over embeddings of bases and
$\M P\M B_m(G)$  be
the category of smooth  principal  $G$-bundles  with  $m$-dimensional
bases  and
smooth $G$-bundle maps $(\varphi,f)$, where the map $f \in$
Mor$\M M_m$.

\bDf\label{Df1.1}
A {\it $G$-gauge-natural bundle\/}  is
a covariant functor $F$ from the category $\M P\M B_m(G)$ to
the category $\M F\M M_m$  satisfying

   i) for each $\pi:\f P \to \f M$ in $\Ob\M P\M B_m(G),
\pi_{\f P}:F\f P \to \f M$
 is  a  fibered  manifold over $\f M$,
\par   ii) for each map $(\varphi,f)$ in $\Mor\M P\M B_m(G),
F\varphi = F(\varphi,f)$
  is  a  fibered  manifold morphism covering $f$,

   iii) for any open subset $\f U \subseteq \f M$, the immersion
$\iota:\pi^{-1}(\f U)\hookrightarrow \f P$  is  transformed
      into the immersion $F\iota:\pi^{-1}_{\f P}(\f U)
\hookrightarrow F\f P$.
\eDf

Let $(\pi:\f P \to \f M) \in\Ob \M P\M B_m(G)$ and $W^r\f P$ be
the space of all $r$-jets $j^r_{(0,e)}\varphi$,  where
$\varphi:\Rn^m\times G \to \f P$ is
in $\Mor \M P\M B_m(G), 0\in \Rn^m$  and $e$
is the unit in $G$. The space $W^r\f P$ is  a
principal fiber bundle over the manifold $\f M$ with
the  structure  group
$W^r_m G$ of  all  $r$-jets $j^r_{(0,e)}\Psi$
of  principal fiber  bundle  isomorphisms
$\Psi:\Rn^m\times G \to \Rn^m\times G$ covering the diffeomorphisms
$\psi:\Rn^m \to \Rn^m$
  such  that  $\psi(0)=0$.  The
group $W^r_mG$ is the semidirect product of $G^r_m
=\inv J^r_0(\Rn^m,\Rn^m)_0$  and $T^r_mG= J^r_0(\Rn^m,G)$
 with respect to the  action
of $G^r_m$  on $T^r_mG$ given by the jet composition, i.e.
$W^r_m G= G^r_m\sem T^r_m G$.
If $(\varphi:\f P \to \bar{\f P}) \in \Mor \M P\M B_m(G)$,  then  we
can  define  the  principal  bundle  morphism
$W^r\varphi:W^r\f P \to W^r\bar{\f P}$
   by   the   jet
composition. The rule transforming any $\f P \in \Ob \M P\M B_m(G)$
into  $W^r\f P \in \Ob \M P\M B_m(W^r_mG)$  and
any $\varphi \in \Mor \M P\M B_m(G)$ into $W^r\varphi \in \Mor
\M P\M B_m(W^r_mG)$ is a $G$-gauge-natural bundle.

The gauge-natural bundle functor $W^r$  plays a  fundamental  role  in
the theory of gauge-natural bundles. We have, \cite{Eck81, KolMicSlo93},

\bTh\label{Th1.2}
Every gauge-natural bundle
is a fiber bundle associated to
the bundle $W^r$  for a certain order $r$.
\eTh

The number $r$ from Theorem \ref{Th1.2} is  called
the {\it order  of  the  gauge-natural
bundle}. So if $F$ is an $r$-order gauge-natural bundle, then
\beq
F\f P = (W^r\f P,S_F), \quad F\varphi = (W^r\varphi,\text{id}_{S_F}),
\eeq
where $S_F$  is a left $W^r_mG$-manifold called the
{\it standard fiber} of $F$.

If $(x^\lambda,z^a)$ is a local fiber coordinate chart on $\f P$
and $(y^i)$ a coordinate chart on
$S_F$, then $(x^\lambda,y^i)$ is the fiber coordinate chart on $F\f P$
 which is said to be {\it adapted}.

Let $F$ be a $G$-gauge-natural bundle of order $s$ and let
$r \le s$ be a minimal number such that the action of
$W^s_mG = G^s_m\sem T^s_mG$ on  $S_F$
can be factorized through the canonical projection
$\pi^s_r:T^s_mG \to T^r_mG.$
Then $r$ is called the {\it gauge-order} of $F$
and we say that $F$ is of order $(s,r)$. We
shall denote by $W_m^{(s,r)}G = G^s_m\sem T^r_mG$ the
Lie group acting on the standard fiber of
an $(s,r)$-order $G$-gauge-natural bundle.
Then there is a one-to-one,  up  to
equivalence, correspondence of smooth left $W^{(s,r)}_mG$-manifolds
and  $G$-gauge-natural
bundles of order $(s,r)$, \cite{Eck81}. So any $(s,r)$-order
$G$-gauge-natural bundle  can  be
represented by its standard fiber with an action
of the group $W^{(s,r)}_mG$.

If $F$ is an $(s,r)$-order $G$-gauge-natural bundle, then $J^kF$
is an $(s+k,r+k)$-order $G$-gauge-natural bundle with the standard fiber
$T^k_mS_F= J^k_0(\Rn^m, S_F)$.

The class of $G$-gauge-natural bundles contains the class
of  natural  bundles in the sense of
\cite{KolMicSlo93, KruJan90, Nij72, Ter78}. Namely, if $F$ is an
$r$-order
natural bundle, then $F$ is the $(r,0)$-order
$G$-gauge-natural bundle  with  trivial
gauge structure.

\smallskip

Let $F$ be a $G$-gauge-natural bundle and $(\varphi,f):\f P
\to \bar{\f P}$
be in the category $\M P\M B_n(G)$.
Let $\sigma $ be a section of $F\f P$. Then we define the section
$\varphi ^*_F\sigma =  F\varphi \circ \sigma \circ f^{-1}$    of
$ F\bar{\f P}$.
Let $H$ be another gauge-natural bundle.

\bDf\label{Df1.3}
A {\it natural differential
operator\/} from $F$
to $H$ is a collection $D = \{D(\f P),\f P\in\Ob\M P\M B_n(G)$\}
of differential operators from $C^\infty(F\f P)$
to $C^\infty(H\f P)$ satisfying
$D(\bar{\f P})\circ \varphi ^*_F  = \varphi ^*_H
\circ D(\f P)$ for each map $(\varphi ,f) \in \Mor\M P\M
B_n(G), \varphi :\f P \to \bar{\f P}$.
\eDf

$D$
is of order $k$ if all $D(\f P)$ are of order $k$.
Let $D$ be a natural differential operator of order $k$ from
$F$  to  $H$.  For  any
$\f P \in\Ob\M P\M B_n(G)$ we have the associated map
$\M D(\f P):J^kF\f P \to H\f P$, over $\f M$, defined  by
$\M D(\f P)(j^k_x\sigma) =  D(\f P)\sigma (x)$  for  all $ x \in \f M$
and  any  section  $\sigma :\f M \to F\f P$.  From  the
naturality of $D$ it follows that $\M D = \{\M D(\f P),\f P
\in \Ob\M P\M B_n(G)\}$ is a natural  transformation of
$J^kF$ to $H$. The following theorem is due to Eck, \cite{Eck81}.

\bTh\label{Th1.4}
Let $F$ and $H$ be $G$-gauge-natural bundles
of  order  $\le (s,r),  s \ge r$.
Then we have a one-to-one correspondence between natural differential
operators of order $k$ from $F$ to $H$ and
$W^{(s+k,r+k)}_nG$-equivariant
 maps from $T^k_m S_F$  to $S_H$.
\eTh

So according to Theorem \ref{Th1.4} a classification of natural
operators between $G$-gauge-natural bundles is equivalent
to the classification of equivariant maps between standard fibers.
Very important tool in classifications of equivariant maps is
the {\em orbit reduction theorem}, \cite{KolMicSlo93, KruJan90}.
Let $p: G\to H$ be a Lie group epimorphism with the kernel $K$,
$M$ be a left $G$-space, $Q$ be a left $H$-space and $\pi: M\to Q$ be
a $p$-equivariant surjective submersion, i.e., $\pi(gx) = p(g) \pi(x)$
for all $x\in M$, $g\in G$. Having $p$, we can consider every  left
$H$-space $N$ as a left $G$-space by $gy=p(g)y$, $g\in G$, $y\in N$.

\bTh\label{Th1.5}
If each $\pi^{-1}(q)$, $q\in Q$ is a $K$-orbit in $M$, then there
is a bijection between the $G$-maps $f: M\to N$ and the $H$-maps
$\varphi : Q\to N$ given by $f=\varphi\com \pi$.
\eTh

\section{Linear connections on vector bundles}
\label{Linear connections on vector bundles}
\setcounter{equation}{0}
\setcounter{theorem}{0}

In what follows let $G=GL(n,\Rn)$ be the group of linear automorphisms
of $\Rn^n$ with coordinates $(a^i_j)$.
Let us consider the category $\M{VB}_{m,n}$
of vector bundles with $m$-dimensional bases,
$n$-dimensional fibers and local fibered linear diffeomorphisms.
Then any vector bundle
$
(p: \f E\to \f M) \in\Ob\M{VB}_{m,n}
$
can be considered as a zero
order $G$-gauge-natural functor $\M{PB}_{m}(G)\to \M{VB}_{m,n}$.

Local linear fiber coordinate
charts on $\f E$ will be denoted by $(x^\lam, y^i)$.
The induced
local bases of sections of
$T\f E$ or
$T^*\f E$
will be denoted by
$(\der_\lam, \der_i)$ or $(d^\lam,d^i)$, respectively.

\smallskip

We define a {\em linear connection\/} on $\f E$ to be a linear splitting
\beq
K : \f E \to J^1\f E\,.\quad
\eeq
Considering the contact morphism $J^1\f E\to
T^*\f M \ten T\f E$ over the identity of $T\f M$, a
linear connection can be regarded as a $T\f E$-valued
1-form
\beq
K : \f E \to T^*\f M \ten T\f E
\eeq
projecting onto the identity of $T\f M$.

The coordinate expression of a linear connection
$K$
is of the type
\beq
K = d^\lam \ten \big(\der_\lam
+ K\col ji\lam \, y^j \, \der_i \big) \,,
\sep{with}  K\col ji \lam \in C^\infty(\f M,\Rn)\,.
\eeq

Linear connections can be regarded as sections of
a (1,1)-order $G$-gauge-natural bundle
$\Lin\f E \to \f M$, \cite{Eck81, KolMicSlo93}.
The standard fiber of the functor $\Lin$ will be denoted by
${R}=\Rn^{n*}\ten \Rn^n\ten\Rn^{m*}$, elements of $R$ will be
said to be {\em formal linear connections},
the induced coordinates
on ${R}$ will be said to be {\em formal symbols} of formal linear
connections and will be denoted by $(K\col ji\lam)$.
The action $\bet:W^{(1,1)}_mG \car {R} \to {R}$
of the group $W^{(1,1)}_mG=G^1_m\sem T^1_mG$
on the standard fiber ${R}$
is given in coordinates by
\beq
(K\col ji\lam) \com \bet=
        a^i_p \,( K\col qp\rho
        \tilde a^q_j \tilde a^\rho_\lam
        -\tilde a^p_{j\lam})\,,
\eeq
where $(a^\lam_\mu,a^i_j,a^i_{j\lam})$ are coordinates on
$W^{(1,1)}_mG$ and $\tilde{\ }$ denotes the inverse element.

\bNt\label{Nt2.1}
Let us note that the action $\bet$ gives, in a natural way,
the action
\beq
\bet^r:W^{(r+1,r+1)}_m G\car T^r_m{R} \to T^r_m{R} \,
\eeq
determined by the jet prolongation of the action
$\bet$.
\eNt

\bRm\label{Rm2.2}
Let us consider the group epimorphism
$\pi^{r+1,r+1}_{r,r}: W^{(r+1,r+1)}_mG\to  W^{(r,r)}_mG$
and its kernel $B^{r+1,r+1}_{r,r} G\byd\Ker \pi^{r+1,r+1}_{r,r} $.
On $B^{r+1,r+1}_{r,r} G$ we have the induced
coordinates $(a^\lam_{\mu_1\dots\mu_{r+1}},
a^i_{j\mu_1\dots\mu_{r+1}})$.
 Then the restriction $\bar\bet^r$ of the action
$\bet^r$ to $B^{r+1,r+1}_{r,r} G$ has the following coordinate
expression
\begin{align}\label{Eq2.1}
(K\col ji{\mu_1} , &  \dots,
        K\col ji{\mu_1,\mu_2\dots\mu_{r+1}}) \com \bar\bet^r
\\
 & =
        (K\col ji{\mu_1},\dots,
        K\col ji{\mu_1,\mu_2\dots\mu_{r}},
        K\col ji{\mu_1,\mu_2\dots\mu_{r+1}}-
        \tilde a^i_{j\mu_1\dots\mu_{r+1}})\,,\nonumber
\end{align}
where $(K\col ji{\mu_1} ,
K\col ji{\mu_1,\mu_2},  \dots,
        K\col ji{\mu_1,\mu_2\dots\mu_{r+1}})$ are
the induced jet coordinates on $T^r_m{R}$.
\eRm

The curvature of a linear connection
${K}$
on
$\f E$
turns out to be the vertical valued 2--form
\beq
R[{K}] = - [{K},{K}]:\f E\to V\f E\ten \bigwedge^2 T^*\f M\,,
\eeq
where $[,]$ is the Froelicher-Nijenhuis bracket. The coordinate
expression is
\begin{align*}
R[{K}] & = R[K]\col ji{\lam\mu}\, y^j\,  \der_i \ten d^\lam \wed d^\mu
\\
& = -2(\der_\lam K \col ji\mu + K\col jp\lam K\col pi\mu)\, y^j\,
        \der_i \ten d^\lam \wed d^\mu \,.
\end{align*}

If we consider the identification $V\f E = \f E \ucar{\f M}\f E$
and linearity of $R[K]$, the curvature $R[K]$ can be considered as
the curvature tensor field
$
R[{K}] :\f M\to\f E^* \ten \f E \ten \bigwedge^2 T^*\f M\,
$
and
\beq
R[{K}]:C^\infty(\Lin\f E) \to
        C^\infty(\f E^*\ten\f E\ten\bigwedge^2T^*\f M)
\eeq
is a natural operator which is of order one,
i.e., we have the associated
$W^{(2,2)}_m G$-equivariant map, called the {\em formal curvature
map}
of formal linear connections,
$$
\M R_L: T^1_m {R} \to \M U
$$
with the coordinate expression
\bEq\label{Eq2.2}
(u_j{}^i{}_{\lam\mu})\com\M R_L= K \col ji{\lam,\mu}
        - K \col ji{\mu,\lam} + K\col jp\mu K\col pi\lam
        - K\col jp\lam K\col pi\mu \,,
\eEq
where $(u_j{}^i{}_{\lam\mu})$ are the induced
coordinates on the standard fiber
$\M U\byd \Rn^{n*}\ten\Rn^n\ten\bigwedge^2 \Rn^{m*}$
of  $\f E^*\ten\f E\ten\bigwedge^2T^*\f M$.


\smallskip

We define a {\em classical connection\/} on $\f M$ to be
a linear symmetric connection on the tangent vector bundle
$p_{\f M}:T\f M\to \f M$ with the coordinate
expression
\beq
\Lam = d^\lam \ten \big(\der_\lam
+ \Lam\col\nu\mu\lam \, \dot x^\nu \, \dt\der_\mu \big) \,,
\quad  \Lam\col \mu \lam \nu\in C^\infty(\f M,\Rn),\quad
\Lam\col \mu \lam \nu = \Lam\col \nu \lam \mu \,.
\eeq

Classical connections can be regarded as sections of a 2nd
order natural bundle
$\Cla\f M \to \f M$, \cite{KolMicSlo93}.
The standard fiber of the functor $\Cla$ will be denoted by
${Q}=\Rn^m\ten S^2\Rn^{m*}$, elements of $Q$ will be said
to be {\em formal classical connections},
the induced
coordinates on ${Q}$ will be said to be {\em formal Christoffel
symbols} of formal classical connections and will be denoted by
$(\Lam\col \mu \lam \nu)$.
The action $\alp:G^2_m \car {Q} \to {Q}$
of the group $G^2_m$ on ${Q}$
is given in coordinates by
\beq
(\Lam\col \mu \lam \nu) \com \alp=
        a^\lam_\rho \,( \Lam\col \sig \rho \tau
        \tilde a^\sig_\mu \tilde a^\tau_\nu
        -\tilde a^\rho_{\mu\nu})\,.
\eeq

\bNt\label{Nt2.3}
Let us note that the action $\alp$ gives, in a natural way,
the action
\beq
\alp^r:G^{r+2}_m \car T^r_m{Q} \to T^r_m{Q} \,
\eeq
determined by the jet prolongation of the action
$\alp$.
\eNt

\bRm\label{Rm2.4}
Let us consider the group epimorphism
$\pi^{r+2}_{r+1}: G^{r+2}_m\to  G^{r+1}_m$
and its kernel $B^{r+2}_{r+1}\byd\Ker \pi^{r+2}_{r+1}$. We have
the induced coordinates $(a^\lam_{\mu_1\dots\mu_{r+2}})$
on $B^{r+2}_{r+1}$. Then the restriction $\bar\alp^r$ of the action
$\alp^r$ to $B^{r+2}_{r+1}$ has the following coordinate
expression
\begin{align}\label{Eq2.3}
(\Lam\col {\mu_1} \lam {\mu_2}, & \dots,
        \Lam\col {\mu_1} \lam {\mu_2,\mu_3\dots\mu_{r+2}}) \com
        \bar\alp^r
\\
 & =
        (\Lam\col {\mu_1} \lam {\mu_2},\dots,
        \Lam\col {\mu_1} \lam {\mu_2,\mu_3\dots\mu_{r+1}},
        \Lam\col {\mu_1} \lam {\mu_2,\mu_3\dots\mu_{r+2}}-
        \tilde a^\lam_{\mu_1\dots\mu_{r+2}})\,,\nonumber
\end{align}
where $(\Lam\col {\mu_1} \lam {\mu_2},
\Lam\col {\mu_1} \lam {\mu_2,\mu_3},  \dots,
        \Lam\col {\mu_1} \lam {\mu_2,\mu_3\dots\mu_{r+2}})$ are
the induced jet coordinates on $T^r_m{Q}$.
\eRm

\smallskip

The curvature tensor of a classical connection is a natural operator
\beq
R[\Lam]:C^\infty(\Cla\f M) \to C^\infty(T^*\f M\ten T\f M\ten
        \bigwedge^2T^*\f M)
\eeq
which is of order one, i.e., we have the associated
$G^3_m$-equivariant map, called the {\em formal curvature
map} of formal classical connections,
\beq
\M R_C:T^1_m{Q} \to S_{T^*\ten T\ten
        \bigwedge^2T^*}
\eeq
with the coordinate expression
\beq
(w_\nu{}^\rho{}_{\lam\mu})\com\M R_C= \Lam\col \nu\rho{\lam,\mu}
        - \Lam\col \nu \rho {\mu,\lam} + \Lam\col \nu\sig \mu
        \Lam\col\sig\rho\lam - \Lam\col \nu\sig \lam
        \Lam\col\sig\rho\mu \,,
\eeq
where $(w_\nu{}^\rho{}_{\lam\mu})$ are the induced
coordinates on the standard fiber $\M W\byd S_{T^*\ten T\ten
\bigwedge^2T^*}
=\Rn^{m*}\ten\Rn^m\ten\bigwedge^2\Rn^{m*}$.

\smallskip

Let us denote by
$
\f E^{p,r}_{q,s}\byd \ten^p \f E
\ten\ten^q \f E^*\ten\ten^rT\f M \ten \ten^s T^*\f M
$
the tensor product over $\f M$ and recall that
$
\f E^{p,r}_{q,s}
$
is a vector bundle which is a $G$-gauge-natural bundle of order
$(1,0)$.

A classical connection $\Lam$ on $\f M$ and a linear
connection ${K}$ on $\f E$ induce the linear tensor product connection
$K^p_q\ten {\Lam}^r_s\byd \ten^p K\ten\ten^q K^*\ten
\ten^r{\Lam} \ten \ten^s{\Lam}^*$ on
$\f E^{p,r}_{q,s}$
\beq
K^p_q\ten {\Lam}^r_s: \f E^{p,r}_{q,s} \to T^*\f M\uten{\f M}T
        \f E^{p,r}_{q,s}
\eeq
which can be considered
as a linear splitting
\beq
K^p_q\ten{\Lam}^r_s:\f E^{p,r}_{q,s}\to J^1\f E^{p,r}_{q,s}\,.
\eeq

Then we define, \cite{Jan03}, the
{\em covariant differential of a section $\Phi:\f M\to \f E^{p,r}_{q,s}$
with respect to
the pair of connections $(K,\Lam)$} as a section of
$\f E^{p,r}_{q,s}\ten T^*\f M$ given by
\beq
\nab^{(K,\Lam)} \Phi = j^1 \Phi - (K^p_q\ten{\Lam}^r_s )\com \Phi\,.
\quad
\eeq

In what follows we set $\nab = \nab^{(K,\Lam)}$ and
$\phi^{i_1\dots i_p\lam_1\dots\lam_r}
        _{j_1 \dots j_q\mu_1\dots\mu_s;\nu}=
\nab_{\nu}\phi^{i_1\dots i_p\lam_1\dots\lam_r}
        _{j_1 \dots j_q\mu_1\dots\mu_s}$.

\smallskip

We have the following relations between the
covariant differentials and the curvatures, \cite{Jan03}.

\bPr\label{Pr2.5}
The curvature
\beq
R[K^p_q\ten {\Lam}^r_s]\byd - [K^p_q\ten {\Lam}^r_s , K^p_q
                \ten {\Lam}^r_s]
        : \f E^{p,r}_{q,s}\to
        \f E^{p,r}_{q,s}  \ten \bigwedge^2T^*\f M
\eeq
is determined by the curvatures $R[K]$ and $R[{\Lam}]$.
\ePr

\bTh\label{Th2.6}
(The generalized Bianchi identity) We have
\beq
R[{K}]\col {j}i{\lam\mu;\nu}+
        R[{K}]\col {j}i{\mu\nu;\lam}+
        R[{K}]\col {j}i{\nu\lam;\mu}=0\,.
\eeq
\eTh

\bTh\label{Th2.7}
Let $\Phi\in C^\infty(\f E^{p,r}_{q,s})$.
Then we have
\beq
\Alt \nab^2 \Phi =-\frac12\, R[\Lam^p_q \ten{K}^r_s]\com \Phi
        \in C^\infty(\f E^{p,r}_{q,s}\ten\bigwedge^2T^*\f M)\,,
\eeq
where $\Alt$ is the antisymmetrization.
\eTh

\bRm\label{Rm2.8}
From the above Theorem \ref{Th2.7} and the
expression of $R[K^p_q\ten\Lam^r_s]$, \cite{Jan03},
it follows, that
$\Alt \nab^2 \Phi$ is a $\f E^{p,r}_{q,s}$-valued 2-form
which is a quadratic polynomial in $R[K], R[\Lam], \Phi$.
Especially, we have
\beq
\Alt \nab^2 R[{K}] :\f M\to\f E^*\ten\f E \ten
        \bigwedge^2T^*\f M\ten\bigwedge^2T^*\f M
        \,,
\eeq
given in coordinates by
\begin{align}\nonumber
\Alt \nab^2 R[{K}] & = -\frac12 \, \big(R[K]_{p}{}^i{}_{\nu_1\nu_2}\,
        R[{K}]_j{}^p{}_{\lam\mu} -R[K]_{j}{}^p{}_{\nu_1\nu_2}\,
        R[{K}]_p{}^i{}_{\lam\mu}
\\
&\quad -
        R[\Lam]_{\lam}{}^\ome{}_{\nu_1\nu_2} \,
        R[{K}]_j{}^i{}_{\ome\mu}
        - R[\Lam]_{\mu}{}^\ome{}_{\nu_1\nu_2} \,
        R[{K}]_j{}^i{}_{\lam\ome}\big) \,\nonumber
\\
& \qquad
\E b^j\ten \E b_i\ten d^\lam\wed d^\mu\ten d^{\nu_1}\wed d^{\nu_2}
        \,. \nonumber
\end{align}
\eRm

\bRm\label{Rm2.9}
Let us note that for classical connections we have the first and the
second Bianchi identities
$$
R[{\Lam}]\col {(\nu}\rho{\lam\mu)} = 0\text{\rm \quad and\quad}
        R[{\Lam}]\col {\nu}\rho{(\lam\mu;\sig)}=0\,,
$$
respectively, where $(\dots)$ denotes the cyclic
permutation. Moreover, we have the antisymmetrization
of the second order covariant
differential of the curvature tensor which
is a quadratic polynomial
of the curvature tensor.
\eRm

\section{The first $k$-th order reduction theorem for linear\\
         and classical connections}
\setcounter{equation}{0}
\setcounter{theorem}{0}

Let us introduce the following notations.

Let ${\M W}_0\f M \byd {\M W}\f M=T^*\f M\ten T
\f M\ten\bigwedge^2T^*\f M$,
${\M W}_i\f M= {\M W}\f M\ten\ten^iT^*\f M$, $i\ge 0$. Let us put
${\M W}^{(k,r)}\f M=
{\M W}_k\f M\ucar{\f M} \dots \ucar{\f M} {\M W}_r\f M$.
Especially, we set ${\M W}^{(r)}\f M \byd  {\M W}^{(0,r)}\f M$.
Then ${\M W}_i\f M$ and ${\M W}^{(k,r)}\f M$
are natural bundles of order one and
the corresponding standard
fibers will be denoted by ${\M W}_i$ and ${\M W}^{(k,r)}$,
where ${\M W}_0\byd {\M W} = \Rn^{m*}\ten \Rn^m\ten\bigwedge^2\Rn^{m*}$,
${\M W}_i= {\M W}\ten\ten^i \Rn^{m*}$, $i\ge 0$, and ${\M W}^{(k,r)}=
{\M W}_k\car\dots \car {\M W}_r$.
Let us denote by $(w\col {\nu}{\rho}{\lam\mu\sig_1\dots \sig_i})$
the coordinates on ${\M W}_i$.

We denote by
\beq
\M R_{{C},i}:T^{i+1}_m {Q}\to  {\M W}_i
\eeq
the $G^{i+3}_m$-equivariant map associated with the
$i$-th covariant differential of the curvature tensors of classical
connections
\beq
\nab^i R[\Lam]:C^\infty(\Cla\f M)\to C^\infty({\M W}_i\f M)\,.
\eeq
The map $\M R_{{C},i}$ is said to be the {\em formal
curvature map of order $i$} of classical connections.

Let $C_{{C},i}\subset {\M W}_i$ be a subset given by
identities of the $i$-th covariant differentials of
the curvature tensors of classical connections,
i.e., by covariant differentials of the Bianchi identities
and the antisymmetrization of the second order covariant differentials,
see Remark \ref{Rm2.9}.
So $C_{C,i}$ is given by the following system of equations
\begin{align}
& w\col {(\nu}{\rho}{\lam\mu)\sig_1\dots \sig_i} =0\,,
\label{Eq3.1}
\\
& w\col {\nu}{\rho}{(\lam\mu\sig_1)\sig_2\dots \sig_i} =0\,,
\label{Eq3.2}
\\
& w\col {\nu}{\rho}{\lam\mu\sig_1\dots[\sig_{j-1}
        \sig_j]\dots \sig_i}
        + \pol({\M W}^{(i-2)})=0\,,
\label{Eq3.3}
\end{align}
where $j=2,\dots, i$ and $[..]$ denotes the
antisymmetrization.

Let us put
$C^{(r)}_{{C}}=
C_{{C},0}\car \dots \car C_{{C},r}$
and denote by
$C^{(k,r)}_{C,r^{k-1}_C}$, $k\le r$, the fiber in
$r^{(k-1)}_C\in C^{(k-1)}_C$ of the canonical
projection $\pr^r_{k-1}: C^{(r)}_C\to C^{(k-1)}_C$. For $r < k$
we put $C^{(k,r)}_{C,r^{k-1}_C}=\emptyset$.
Let us note that there is an affine structure on the fibres of
the projection
$\pr^r_{r-1}: C^{(r)}_C\to C^{(r-1)}_C$, \cite{KolMicSlo93}.
Really, $C^{(r)}_C$ is a subbundle in $C^{(r-1)}_C\car {\M W}_r$
given by the solution (for $i=r$) of the system of
nonhomogeneous equations
(\ref{Eq3.1}) -- (\ref{Eq3.3}).

Then we put
\bEq\label{Eq3.4}
\M R^{(k,r)}_C\byd (\M R_{C,k},\dots,\M R_{C,r}):
T^{r+1}_m {Q}\to  {\M W}^{(k,r)}\,,\qquad
\M R^{(r)}_C\byd \M R^{(0,r)}_C\,,
\eEq
which has values,
for any $j^{r+1}_0 \gam\in T^{r+1}_m Q$,
in $C^{(k,r)}_{C,\M R^{(k-1)}(j^{k}_0\gam)}$.
In \cite{KolMicSlo93} it was proved that $C^{(r)}_C$
is a submanifold
in ${\M W}^{(r)}$ and the restriction of $\M R^{(r)}_C$
to $C^{(r)}_C$ is a surjective submersion.
Then we can consider the fiber product
$T^{k}_m Q\car_{C^{(k-1)}_C} C^{(r)}_C$ which will be  denoted
by $T^{k}_m Q\car C^{(k,r)}_C$. In \cite{Jan04a} it was proved
that the mapping
\beq
(\pi^{r+1}_k, \M{R}^{(k,r)}_C): T^{r+1}_m Q
        \to T^k_m Q\car C^{(k,r)}_C
\eeq
is a surjective submersion.

\smallskip

Similarly let  ${\M U}_0\f E\byd{\M U}\f E =
\f E^*\ten\f E\ten\bigwedge^2T^*\f M$,
${\M U}_i\f E= {\M U}\f E\ten\ten^i T^*\f M$, $i\ge 0$,
${\M U}^{(k,r)}\f E=
 {\M U}_k\f E\ucar{\f M}\dots\ucar{\f M}{\M U}_r\f E$. Especially,
${\M U}^{(r)}\f E  \byd {\M U}^{(0,r)}\f E$.
Then ${\M U}_i\f E$ and ${\M U}^{(k,r)}\f E$
are $G$-gauge-natural bundles of order $(1,0)$ and the corresponding
standard fibers will be denoted by ${\M U}_i$ and ${\M U}^{(k,r)}$,
where ${\M U}_0\byd {\M U} = \Rn^{n*}\ten \Rn^n\ten\bigwedge^2\Rn^{m*}$,
${\M U}_i= {\M U}\ten\ten^i\Rn^{m*}$, $i\ge 0$, and ${\M U}^{(k,r)}=
{\M U}_k\car \dots \car {\M U}_r$.
Let us denote by $(u\col {j}{i}{\lam\mu\sig_1\dots \sig_i})$
the coordinates on ${\M U}_i$.

We denote by
\beq
\M R_{{L},i}:T^{i-1}_m Q\car T^{i+1}_m {R}\to  {\M U}_i
\eeq
the ${\M W}^{(i+2,i+2)}_mG$-equivariant map associated with the
$i$-th covariant differential of the curvature tensors of linear
connections
\beq
\nab^i R[K]: C^\infty(\Cla\f M\ucar{\f M}\Lin\f E)
        \to C^\infty({\M U}_i\f E)\,.
\eeq
The map $\M R_{{L},i}$ is said to be the {\em formal
curvature map of order $i$} of general linear
connections.

Let $C_{{L},i}\subset {\M U}_i$ be a subset given by
identities of the $i$-th covariant differentials
of the curvature tensors of linear connections,
i.e., by covariant differentials of the Bianchi identity
and the antisymmetrization of the second order covariant differentials,
see Theorem \ref{Th2.6} and Remark \ref{Rm2.8}.
So $C_{L,i}$ is given by the following system of equations
\begin{align}\label{Eq3.5}
u\col ji{(\lam\mu\sig_1)\rho_2\dots \sig_i} & = 0 \,,
\\
u\col ji{\lam\mu\sig_1 \dots[\sig_{j-1}\sig_j] \dots \sig_i}
        + \pol(C^{(i-2)}_C\car \M U^{(i-2)}) & = 0 \,,\label{Eq3.6}
\end{align}
$j=2,\dots,i$, where
$\pol(C^{(i-2)}_C\car\M U^{(i-2)})$ are some polynomials on
$C^{(i-2)}_C\car\M U^{(i-2)}$.

Let us put
$C^{(r)}_{{L}}=
 C_{{L},0}\car \dots \car C_{{L},r}$
and denote by
$C^{(k,r)}_{L,r^{(k-1)}_L}$, $k\le r$, the fiber in
$r^{(k-1)}_L\in C^{(k-1)}_L$ of the canonical
projection $\pr^r_{k-1}: C^{(r)}_L\to C^{(k-1)}_L$. For $r < k$
we put $C^{(k,r)}_{L,r^{(k-1)}_L}=\emptyset$.
Let us note that there is an affine structure on
the projection
$\pr^r_{r-1}: C^{(r)}_L\to C^{(r-1)}_L$, \cite{Jan04}.
Really, $C^{(r)}_L$ is a subbundle in $C^{(r-1)}_L \car \M U_r$
given as the solution (for $i=r$) of the system of
nonhomogeneous equations (\ref{Eq3.5}) -- (\ref{Eq3.6}).

Then we set
\beq
\M R^{(k,r)}_{{L}}\byd (\M R_{{L},k},\dots,\M R_{{L},r}):
T^{r-1}_m {Q}\car T^{r+1}_m {R}\to {\M U}^{(k,r)}\,,\quad
\M R^{(r)}_L \byd \M R^{(0,r)}_L\,,
\eeq
which has values, for any $(j^{r-1}_0\lam, j^{r+1}_0\gam)\in
T^{r-1}_m Q\car T^{r+1}_m R$,
in $C^{(k,r)}_{{L},\M R^{(k-1)}_L(j^{k-2}_0\lam,j^{k}_0\gam)}$.

In \cite{Jan04} it was proved that
$C^{(s)}_C\car C^{(r)}_L$, $s\ge r-2$, $r\ge 0$,
is a submanifold of
${\M W}^{(s)}\car {\M U}^{(r)}$
and the restriction
\beq
(\M R^{(s)}_C,\M R^{(r)}_L):
T^{s+1}_m {Q}\car T^{r+1}_m {R}\to C^{(s)}_C\car C^{(r)}_L
\eeq
is a surjective submersion.
Then we can consider the fiber product
\beq
(T^{k_1}_m Q\car T^{k_2}_m R) \ucar{C^{(k_1-1)}_C\car C^{(k_2-1)}_L}
        (C^{(s)}_C\car C^{(r)}_L)\,,
\eeq
$k_1\ge k_2-2$, and  denote it
by $T^{k_1}_m Q\car T^{k_2}_m R \car C^{(k_1,s)}_C\car C^{(k_2,r)}_L$.

\smallskip

Now we shall prove the technical

\bLm\label{Lm3.1}
If $s\ge r-2$, $k_1\ge k_2-2$, $s+1\ge k_1$,
$r+1\ge k_2$,  then
the restricted map
\beq
(\pi^{s+1}_{k_1}\car \pi^{r+1}_{k_2},\M R^{(k_1,s)}_C,\M R^{(k_2,r)}_L):
T^{s+1}_m {Q}\car T^{r+1}_m {R}\to T^{k_1}_m Q\car T^{k_2}_m R
        \car  C^{(k_1,s)}_C\car C^{(k_2,r)}_L
\eeq
is a surjective submersion.
\eLm

\bPf
In \cite{Jan04a} it was proved that
\beq
 (\pi^{s+1}_{k_1}, \M R^{(k_1,s)}_C):T^{s+1}_m Q \to
        T^{k_1}_m Q\car  C^{(k_1,s)}_C
\eeq
is a surjective submersion.
The mapping of Lemma \ref{Lm3.1} is then a surjective submersion if and
only if the mapping $(\pi^{r+1}_{k_2}, \M R^{(k_2,r)}(j^{s+1}_0\lam,-))
:T^{r+1}_m R \to T^{k_2}_m R\car C^{(k_2,r)}_L$ is a surjective
submersion for any $j^{s+1}_0\lam\in T^{s+1}_m Q$.
Let us assume $i=k_2,\dots, r$.
By \cite{Jan04} the mapping $\M R^{(i)}_L(j^{s+1}_0\lam, -):
T^{i+1}_m R \to C^{(i)}_L$ is a surjective submersion and we have the
commutative diagram
\beq
\begin{CD}
T^{i+1}_m R
        @> \M R^{(i)}_L(j^{s+1}_0\lam,-)>>
        C^{(i)}_L
\\
        @V{\pi ^{i+1}_i}VV
        @VV{\pr^i_{i-1}}V
\\
        T^{i}_m R
        @>\M R^{(i-1)}_L(j^{s+1}_0\lam,-)>>
        C^{(i-1)}_L
\end{CD}
\eeq
All morphisms in the above diagram are surjective submersions
which implies that the mapping
$(\pi^{i+1}_i,\M R^{(i)}_L(j^{s+1}_0\lam,-)):
T^{i+1}_m {R}\to  T^i_m R
        \ucar{C^{(i-1)}_L}
         C^{(i)}_L
$
is a surjection
over $\M R^{(i-1)}_L(j^{s+1}_0\lam,-)$ given by
$(\pi^{i+1}_i,\M R_{L,i}(j^{s+1}_0\lam,-))$.
But the mapping $\M R_{L,i}(j^{s+1}_0\lam,-)$ is affine
morphisms over $\M R^{(i-1)}_L(j^{s+1}_0\lam,-)$
(with respect to the
affine structures on
$\pi^{i+1}_i:T^{i+1}_m R
\to  T^i_m R$ and
$\pr^i_{i-1}: C^{(i)}_L\to  C^{(i-1)}_L$)
which has a constant rank.
So the surjective morphism
$
(\pi^{i+1}_i,\M R_{L,i}(j^{s+1}_0\lam,-))
$
has a constant rank and hence is a submersion.
$
(\pi^{r+1}_{k_2},\M R^{(k_2,r)}_L(j^{s+1}_0\lam,-))
$
is then a composition of surjective submersions
\begin{align*}
(\pi^{k_2+1}_{k_2},
\M R_{L,k_2}(j^{s+1}_0\lam,-), \id_{C^{(k_2+1,r)}_L}) & \com  \dots
\\
\dots\com (\pi^{r}_{r-1}, \M R_{L,r-1}(j^{s+1}_0\lam,-),
\id_{C^{(r,r)}_L})
 & \com (\pi^{r+1}_r,\M R_{L,r}(j^{s+1}_0\lam,-))\,.
\end{align*}
\vglue-1.3\baselineskip

\ePf

Let $F$ be a $G$-gauge-natural bundle of order $k$, i.e.,
$S_F$ is a $W^{(k,k)}_mG$-manifold.

\bTh\label{Th3.2}
Let $s\ge r-2$, $r+1, s+2\ge k\ge 1$.
For every $W^{(s+2,r+1)}_m G$-equivariant map
$$
f:T^s_m Q\car T^r_m R\to S_F
$$
there exists a unique
$W^{(k,k)}_m G$-equivariant map
$$
g:T^{k-2}_m Q\car T^{k-1}_m R\car C^{(k-2,s-1)}_C \car
        C^{(k-1,r-1)}_L \to S_F
$$
satisfying
\beq
f=g\com (\pi^{s}_{k-2}\car \pi^{r}_{k-1},
\M R^{(k-2,s-1)}_C, \M R^{(k-1,r-1)}_L)\,.
\eeq
\eTh

\bPf
Let us consider the space
\beq
S_{C,s} \byd \Rn^m\ten S^s\Rn^{m*}\quad\text{or}\quad
        S_{L,r} \byd \Rn^{n*}\ten\Rn^{n}\ten S^r \Rn^{m*}\,
\eeq
with coordinates $(s^\lam{}_{\mu_1\mu_2\dots\mu_{s}})$
or $(s\col ji{\mu_1\dots\mu_{r}})$, respectively. Let us consider
the action of $G^s_m$ on $S_{C,s}$ and the action of
$W^{(r,r)}_m G$ on $S_{L,r}$ given by
\bEq\label{Eq3.7}
\bar s^\lam{}_{\mu_1\mu_2\dots\mu_{s}}=
         s^\lam{}_{\mu_1\mu_2\dots\mu_{s}}
        - \tilde a^\lam_{\mu_1\dots\mu_{s}}\,,\quad
\bar s\col ji{\mu_1\dots\mu_{r}}=
         s\col ji{\mu_1\dots\mu_{r}}
        - \tilde a^i_{j\mu_1\dots\mu_{r}}\,.
\eEq
From (\ref{Eq2.1}), (\ref{Eq2.3})
and (\ref{Eq3.7}) it is easy to see that
the symmetrization maps
\beq
\sig_{C,s}:T^{s}_m {Q} \to  S_{C,s+2}\,,\quad
\sig_{L,r}:T^{r}_m {R} \to S_{L,r+1}
\eeq
given by
\beq
(s^\lam{}_{\mu_1\mu_2\dots\mu_{r+1}})\com
\sig_{C,s} = \Lam\col {(\mu_1}\lam{\mu_2,\mu_3\dots\mu_{s+2})}\,,
\quad
(s\col ji{\mu_1\dots\mu_{r+1}})\com
\sig_{L,r} = K\col ji{(\mu_1,\mu_2\dots\mu_{r+1})}\,
\eeq
are equivariant.

We have the $G^{s+2}_m$-equivariant
map
\begin{align*}
\varphi_{C,s} & \byd(\sig_{C,s},\pi^{s}_{s-1},
        \M R_{C,s-1})
\\
      & :  T^{s}_m {Q}  \to
        S_{C,s+2}  \car  T^{s-1}_m {Q}\car  {\M W}_{s-1}\,.
\end{align*}
On the other hand we define the $G^{s+2}_m$-equivariant map
\begin{align*}
\psi_{C,s} & : S_{C,s+2} \car
        T^{s-1}_m {Q}\car
        {\M W}_{s-1}\to
        T^{s}_m {Q}
        \,
\end{align*}
over the identity of $T^{s-1}_m {Q}$
by the following coordinate expression
\begin{align}
\Lam\col \mu\lam{\nu,\rho_1\dots\rho_{s}} & =
        s^\lam{}_{\mu\nu\rho_1\dots\rho_{s}}+\lin(
        w_\mu{}^\lam{}_{\nu\rho_1\dots\rho_{s}}
        -\pol(T^{s-1}_m Q))\,,\label{Eq3.8}
\end{align}
where $\lin$ denotes the linear combination with real coefficients
which arises in the following way.
We recall that  $\M R_{C,s-1}$ gives the coordinate
expression
\begin{align}
\Lam\col \mu\lam{\nu,\rho_1\dots\rho_{s}} -
        \Lam\col \mu\lam{\rho_1,\nu\rho_2\dots\rho_{s}}
        & = w_\mu{}^\lam{}_{\nu\rho_1\dots\rho_{s}} -
        \pol(T^{s-1}_m {Q}))\,.\label{Eq3.9}
\end{align}
We can write
\begin{align*}
\Lam\col \mu\lam{\nu,\rho_1\dots\rho_{s}} & =
        s^\lam{}_{\mu\nu\rho_1\dots\rho_{s}}+
        (\Lam\col \mu\lam{\nu,\rho_1\dots\rho_{s}} -
        \Lam\col {(\mu}\lam{\nu,\rho_1\dots\rho_{s})})\,.
\end{align*}
Then the term in brackets can be written as a linear
combination of terms of the type
\beq
\Lam\col \mu\lam{\nu,\rho_i\rho_1\dots\rho_{i-1}\rho_{i+1}\dots
        \rho_{s}} -
        \Lam\col \mu\lam{\rho_i,\nu\rho_1\dots\rho_{i-1}\rho_{i+1}\dots
        \rho_{s}}\,,
\eeq
$i=1,\dots,s$, and from (\ref{Eq3.9}) we get (\ref{Eq3.8}).

Moreover,
\beq
\psi_{C,s}\com\varphi_{C,s}=\id_{T^{s}_m {Q}}\,.
\eeq

Similarly we have the $W^{(r+1,r+1)}_m G$-equivariant
map
\begin{align*}
\varphi_{L,r} & \byd(\sig_{L,r},\id_{T^{r-2}_m Q}\car\pi^r_{r-1},
        \M R_{{L,r-1}})
\\
      & :  T^{r-2}_m {Q} \car T^{r}_m {R}  \to
        S_{L,r+1} \car
        T^{r-2}_m {Q}\car T^{r-1}_m {R}  \car
        {\M U}_{r-1}\,
\end{align*}
and we define the $W^{(r+1,r+1)}_m G$-equivariant
map
\begin{align*}
\psi_{L,r} & : S_{L,r+1} \car
        T^{r-2}_m {Q}\car T^{r-1}_m {R}  \car
        {\M U}_{r-1} \to
        T^{r-2}_m {Q} \car T^{r}_m {R}
        \,
\end{align*}
over the identity of $T^{r-2}_m {Q}\car T^{r-1}_m {R}$
by the following coordinate expression
\begin{align}
 K\col ji{\lam,\rho_1\dots\rho_{r}} & =
        s\col ji{\lam\rho_1\dots\rho_{r}} + \lin(
        u\col ji{\lam\rho_1\dots\rho_{r}}
       -\pol(T^{r-2}_m Q\car T^{r-1}_m R))\,,\label{Eq3.10}
\end{align}
where $\lin$ denotes the linear combination with real coefficients
which arises in the following way.
We recall that  $\M R_{{L},r-1}$ gives the coordinate
expression
\begin{align}
      K\col ji{\lam,\rho_1\dots\rho_{r}} -
         K\col ji{\rho_1,\lam\rho_2\dots\rho_{r}}
= u\col ji{\lam\rho_1\dots\rho_{r}}
  - \pol(T^{r-2}_m {Q}\car T^{r-1}_m {R}))\,.\label{Eq2.11}
\end{align}
We can write
\begin{align*}
 K\col ji{\lam,\rho_1\dots\rho_{r}} & =
        s\col ji{\lam\rho_1\dots\rho_{r}} +
        ( K\col ji{\lam,\rho_1\dots\rho_{r}}-
         K\col ji{(\lam,\rho_1\dots\rho_{r})})\,.
\end{align*}
Then the term in brackets can be written as a linear
combination of terms of the type
\begin{align*}
&
         K\col ji{\lam,\rho_i\rho_1\dots\rho_{i-1}
                \rho_{i+1}\dots\rho_{r}} -
         K\col ji{\rho_i,\lam\rho_1\dots\rho_{i-1}
                \rho_{i+1}\dots\rho_{r}}\,,
\end{align*}
$i=1,\dots,r$, and from(\ref{Eq2.11}) we get (\ref{Eq3.10}).

Moreover,
\beq
\psi_{L,r}\com\varphi_{L,r}=\id_{T^{r-2}_m {Q}
\car T^{r}_m {R} }\,.
\eeq

Now we have to distinguish three possibilities.

\smallskip

A) Let $s=r-1$. We have the same orders of groups $G^{r+1}_m$
and $W^{(r+1,r+1)}_m G$ acting on $T^{r-1}_m Q$ and $T^r_m R$.

Let us denote by
\beq
A^r\byd
        T^{r-2}_m {Q}\car T^{r-1}_m {R}  \car
        {\M W}_{r-2}\car  {\M U}_{r-1}\,.
\eeq
Then the map $f\com(\psi_{C,r-1},\psi_{L,r}):
S_{C,r+1}\car  S_{L,r+1} \car A^r\to S_F$ satisfies
the conditions of
the orbit reduction Theorem \ref{Th1.5} for the group epimorphism
$\pi^{r+1,r+1}_{r,r}:
W^{(r+1,r+1)}_m G\to W^{(r,r)}_m G$ and the surjective submersion
$\proj_3:S_{C,r+1}\car  S_{L,r+1} \car A^r\to A^r$. Indeed, the
space
$S_{C,r+1}\car
S_{L,r+1}$ is a $B^{r+1,r+1}_{r,r} G$-orbit. Moreover,
(\ref{Eq3.7}) implies that the action
of $B^{r+1,r+1}_{r,r} G$ on
$S_{C,r+1}\car S_{L,r+1}$ is simply transitive.
Hence there exists a unique $W^{(r,r)}_m G$-equivariant
map
\beq
g_r: A^r=T^{r-2}_m {Q}\car T^{r-1}_m {R}
        \car {\M W}_{r-2}
        \car {\M U}_{r-1} \to S_F
\eeq
such that the following diagram
\begin{equation*}\begin{CD}
S_{C,r+1} \car  S_{L,r+1} \car A^r
        @>{(\psi_{C,r-1},\psi_{L,r})}>> T^{r-1}_m {Q} \car T^{r}_m {R}
        @>f>>S_F
\\
@V{\proj_3}VV  @V{(\pi^{r-1}_{r-2}\car\pi^r_{r-1},\M R_{C,r-2},
        \M R_{{{L}},r-1})}VV @V\id_{S_F}VV
\\
A^r @>{\id_{A^r}}>> A^r @>g_r>> S_F
\end{CD}\end{equation*}
commutes. So
$
f\com(\psi_{C,r-1},\psi_{L,r}) = g_r\com \proj_3
$
and if we compose both sides with $(\varphi_{C,r-1},\varphi_{L,r})$, by
considering
$\proj_3\com(\varphi_{C,r-1},\varphi_{L,r})
        = (\pi^{r-1}_{r-2}\car\pi^r_{r-1},\M R_{C,r-2},
        \M R_{{L},r-1})$, we obtain
\beq
f=g_r \com (\pi^{r-1}_{r-2}\car\pi^r_{r-1},\M R_{C,r-2},
        \M R_{{L},r-1})\,.
\eeq

In the second step we consider the same construction for the
map $g_r$ and obtain the commutative diagram
\begin{equation*}\begin{CD}
S_{C,r} \car  S_{L,r} \car A^{r-1}
        \car {\M W}_{r-2}\car {\M U}_{r-1}
        @>{(\psi_{C,r-2},\psi_{L,r-1},\id_{{\M W}_{r-2}\car {\M U}_{r-1}})}>>
        A^r  @>g_r>> S_F
\\
@V{\proj_{3,4,5}}VV  @V{(\pi^{r-2}_{r-3}\car\pi^{r-1}_{r-2},
        \M R_{C,r-3},
        \M R_{{{L}},r-2},\id_{{\M W}_{r-2}\car {\M U}_{r-1}})}VV @V\id_{S_F}VV
\\
A^{r-1}\car {\M W}_{r-2}\car
        {\M U}_{r-1} @>{\id_{A^{r-1}\car {\M W}_{r-2}\car {\M U}_{r-1}}}>> A^{r-1}\car
        {\M W}_{r-2}\car
        {\M U}_{r-1}  @>g_{r-1}>> S_F
\end{CD}\end{equation*}
So that there exists a unique $W^{(r-1,r-1)}_m G$-equivariant
map
$
g_{r-1}: A^{r-1}\car {\M W}_{r-2}
        \car {\M U}_{r-1} \to S_F
$
such that
$g_r= g_{r-1}\com
        (\pi^{r-2}_{r-3}\car\pi^{r-1}_{r-2},\M R_{C,r-3},
        \M R_{{L},r-2}, \id_{{\M W}_{r-2}\car {\M U}_{r-1}} )$,
i.e.
\beq
f=g_{r-1}\com (\pi^{r-1}_{r-3}\car\pi^{r}_{r-2},\M R_{C,r-3},
        \M R_{C,r-2}, \M R_{{L},r-2}, \M R_{{L},r-1})\,.
\eeq

Proceeding in this way we get in the last step
a unique $W^{(k,k)}_m G$-equivariant
map
$$
g_k:T^{k-2}_m Q\car T^{k-1}_m {R} \car {\M W}^{(k-2,r-2)}\car
{\M U}^{(k-1,r-1)}
\to S_F
$$
such that
\beq
f=g_{k}\com (\pi^{r-1}_{k-2}\car \pi^{r}_{k-1},
\M R^{(k-2,r-2)}_{C},\M R^{(k-1,r-1}_{{L}})\,.
\eeq

\smallskip

B) Let $s=r-2$. We have the action of the group $G^{r}_m$
on $T^{r-2}_m Q$
and the action of the group $W^{(r+1,r+1)}_m G$ on $T^r_m R$.

Then the map $f\com(\id_{T^{r-2}_m Q},\psi_{L,r})
:S_{L,r+1} \car T^{r-2}_m Q\car T^{r-1}_m R \car {\M U}_{r-1}\to S_F$
satisfies the conditions of the orbit reduction theorem
\ref{Th1.5} for the group epimorphism
$\pi^{r+1,r+1}_{r,r}:
W^{(r+1,r+1)}_m G\to W^{(r,r)}_m G$ and the surjective submersion
$\pr_{2,3,4}:S_{L,r+1} \car T^{r-2}_m Q\car T^{r-1}_m R \car {\M U}_{r-1}
\to T^{r-2}_m Q\car T^{r-1}_m R \car {\M U}_{r-1}$.
Indeed, the space
$S_{L,r+1}$ is a $B^{r+1,r+1}_{r,r} G$-orbit.
Let us note that the action
of $B^{r+1,r+1}_{r,r} G$ on
$S_{L,r+1}$ is transitive, but not simple transitive.
Hence there exists a unique $W^{(r,r)}_m G$-equivariant
map
$
g_r: T^{r-2}_m {Q}\car T^{r-1}_m {R}
        \car {\M U}_{r-1} \to S_F
$
such that the following diagram
\begin{equation*}
\begin{CD}
        S_{L,r+1} \car T^{r-2}_m Q\car T^{r-1}_m R \car {\M U}_{r-1}
        @>{(\id_{T^{r-2}_m Q},\psi_{L,r})}>>
        T^{r-2}_m {Q} \car T^{r}_m {R}
        @>f>> S_F
\\
        @V{\pr_{2,3,4}}VV
        @V{(\id_{T^{r-2}_m Q}\car \pi^{r}_{r-1},\M R_{L,r-1})}VV
        @V{\id_{S_F}}VV
\\
        T^{r-2}_m Q\car T^{r-1}_m R \car {\M U}_{r-1}
        @>{\id_{T^{r-2}_m Q\car T^{r-1}_m R \car {\M U}_{r-1}}}>>
        T^{r-2}_m Q\car T^{r-1}_m R \car {\M U}_{r-1}
        @>g_r>> S_F
\end{CD}
\end{equation*}
commutes.  So
$
f\com(\id_{T^{r-2}_m Q},\psi_{L,r}) = g_r\com \proj_{2,3,4}
$
and if we compose both sides with
$(\id_{T^{r-2}_m Q},\varphi_{L,r})$, by
considering
$\proj_{2,3,4}\com(\id_{T^{r-2}_m Q},\varphi_{L,r})
        = (\id_{T^{r-2}_m Q}\car\pi^r_{r-1},
        \M R_{{L},r-1})$, we obtain
\beq
f=g_r \com (\id_{T^{r-2}_m Q}\car\pi^r_{r-1},
        \M R_{{L},r-1})\,.
\eeq
Further we proceed as in the second step in A) and we get
a unique
$W^{(k,k)}_m G$-equivariant map
$$
g_k:T^{k-2}_m Q\car T^{k-1}_m {R} \car {\M W}^{(k-2,r-3)}\car
{\M U}^{(k-1,r-1)}
\to S_F
$$
such that
\beq
f=g_k\com (\pi^{r-2}_{k-2}\car \pi^r_{k-1},\M R^{(k-2,r-3)}_{C},
\M R^{(k-1,r-1)}_{{L}})\,.
\eeq

C) Let $s>r-1$. We have the action of the group
$W^{(s+2,r+1)}_m G$ on $T^{s}_m Q\car T^r_m R$.

By \cite{Jan04a} there exists a $W^{(r+1,r+1)}_m G$-equivariant
mapping
$$
g_{r+1}:T^{r-1}_m Q\car T^{r}_m {R} \car {\M W}^{(r-1,s-1)}
\to S_F
$$
such that
$$
  f= g_{r+1}\com (\pi^s_{r-1}\car \id_{T^r_m R}, \M R^{(r-1,s-1)}_C)\,.
$$
$g_{r+1}$ is then the mapping satisfying the condition A), i.e.
there is a unique
$W^{(k,k)}_m G$-equivariant map
$$
g_k:T^{k-2}_m Q\car T^{k-1}_m {R} \car {\M W}^{(k-2,r-2)}\car
{\M U}^{(k-1,r-1)}
\to S_F
$$
such that
\beq
g_{r+1}=g_k\com (\pi^{r-2}_{k-2}\car \pi^r_{k-1},\M R^{(k-2,r-2)}_{C},
\M R^{(k-1,r-1)}_{{L}})\,,
\eeq
i.e.,
\beq
f=g_k\com (\pi^{s}_{k-2}\car \pi^r_{k-1},\M R^{(k-2,s-1)}_{C},
\M R^{(k-1,r-1)}_{{L}})\,,
\eeq

\smallskip
Summarizing all cases we have
\beq
f=g_k\com (\pi^{s}_{k-2}\car \pi^r_{k-1},\M R^{(k-2,s-1)}_{C},
\M R^{(k-1,r-1)}_{{L}})\,
\eeq
for any $s\ge r-2$ and the restriction of $g_k$ to
$T^{r-2}_m Q\car T^{r-1}_m R\car C^{(k-2,s-1)}_{C}\car
C^{(k-1,r-1)}_{{L}}$
is uniquely determined map $g$ we wished to find.
\ePf

In the above Theorem \ref{Th3.2} we have found
a map $g$ which factorizes $f$, but we did not prove, that
$$
(\pi^s_{k-2}\car \pi^r_{k-1}, \M R^{(k-2,s-1)}_{C},
\M R^{(k-1,r-1)}_{{L}}):T^{s}_m {Q}\car
T^{r}_m {R} \to T^{k-2}_m Q\car T^{k-1}_m R\car
C^{(k-2,s-1)}_{C}\car C^{(k-1,r-1)}_{{L}}
$$
satisfy the orbit conditions, namely we did not prove that
$$
(\pi^s_{k-2}\car \pi^r_{k-1}, \M R^{(k-2,s-1)}_{C},
\M R^{(k-1,r-1)}_{{L}})^{-1}(j^{k-2}_0\lam, j^{k-1}_0 \gam,
r^{(k-2,s-1)}_C,r^{(k-1,r-1)}_{{L}})
$$
is a $B^{s+2,r+1}_{k,k}G$-orbit for any $(j^{k-2}_0\lam,
j^{k-1}_0 \gam,r^{(k-2,s-1)}_C,r^{(k-1,r-1)}_{{L}})
\in T^{k-2}_m Q\car T^{k-1}_m R\car
C^{(k-2,s-1)}_{C}\car C^{(k-1,r-1)}_{{L}}$. Now we shall
prove it.

\bLm\label{Lm3.3}
If $(j^s_0 \lam,j^r_0 {\gam}),(j^s_0 \ac{\lam},j^r_0\ac{{\gam}}) \in
T^s_m {Q}\car T^{r}_m {R} $
satisfy
\beq
(\pi^s_{k-2}\car \pi^r_{k-1}, \M R^{(k-2,s-1)}_{C},
\M R^{(k-1,r-1)}_{{L}})(j^s_0 \lam,j^r_0{\gam})=
(\pi^s_{k-2}\car \pi^r_{k-1}, \M R^{(k-2,s-1)}_{C},
\M R^{(k-1,r-1)}_{{L}})(j^s_0\ac{\lam},j^r_0\ac{{\gam}})\,,
\eeq
then there is an element
$h\in B^{s+2,r+1}_{k,k}G$ such that
$h\,.\, (j^s_0 \ac{\lam},j^r_0 \ac{\gam}) = (j^s_0 \lam,j^r_0{\gam})$.
\eLm

\bPf
Consider the orbit set
$(T^s_m {Q}\car T^{r}_m {R} )/B^{s+2,r+1}_{k,k}G$.
This is a $W^{(k,k)}_m G$-set. Clearly the factor projection
\beq
p: T^s_m {Q}\car T^{r}_m {R}
\to (T^s_m {Q}\car T^{r}_m {R} )/ B^{s+2,r+1}_{k,k}G
\eeq
is a $W^{(s+2,r+1)}_m G$-map. By Theorem
\ref{Th3.2} there is a $W^{(k,k)}_m G$-equivariant
map
\beq
g: T^{k-2}_m Q\car T^{k-1}_m R\car
C^{(k-2,s-1)}_{C}\car C^{(k-1,r-1)}_{{L}}   \to
(T^s_m {Q}\car T^{r}_m {R} )/B^{s+2,r+1}_{k,k}G
\eeq
satisfying $p=g\com (\pi^s_{k-2}\car \pi^r_{k-1},
\M R^{(k-2,s-1)}_{C},
\M R^{(k-1,r-1)}_{{L}})$.
If
\begin{multline}
(\pi^s_{k-2}\car \pi^r_{k-1},
\M R^{(k-2,s-1)}_{C},
\M R^{(k-1,r-1)}_{{L}})(j^s_0 \lam,j^r_0  {\gam})
\nonumber
\\
=
(\pi^s_{k-2}\car \pi^r_{k-1},
\M R^{(k-2,s-1)}_{C},
\M R^{(k-1,r-1)}_{{L}})(j^s_0\ac{\lam},j^r_0\ac{{\gam}})
=(j^{k-2}_0 \lam, j^{k-1}_0 \gam, r^{(k-2,s-1)}_C,
r^{(k-1,r-1)}_{{L}})\,,\nonumber
\end{multline}
then
$$
p(j^s_0 \lam,j^r_0 {\gam})=
g(j^{k-2}_0 \lam, j^{k-1}_0 \gam,r^{(k-2,s-1)}_C,r^{(k-1,r-1)}_{{L}})
=p(j^s_0 \ac{\lam},j^r_0 \ac{{\gam}})\,,
$$
i.e. $(j^s_0 \lam,j^r_0 {\gam}), (j^s_0 \ac{\lam},j^r_0 \ac{{\gam}})$
are in the same $B^{s+2,r+1}_{k,k}G$-orbit, which proves
Lemma \ref{Lm3.3}.
\ePf

The space
$T^{k-2}_m Q\car T^{k-1}_m R\car C^{(k-2,s-1)}_{C}\car
C^{(k-1,r-1)}_{{L}}$
is a left $W^{(k,k)}_M G$-space corresponding to the $G$-gauge-natural
bundle
$J^{k-2}\Cla \f M\ucar{\f M} J^{k-1} \Lin \f E\ucar{\f M}
C^{(k-2,s-1)}_C\f M\ucar{\f M} C^{(k-1,r-1)}_L\f E.$
Setting
$
\nab^{(k,s)}=(\nab^k,\dots,\nab^s)
$, then, as a direct consequence of Theorem \ref{Th3.2}, we obtain
the {\em $k$-th order first reduction theorem
for linear and classical connections}.

\bTh\label{Th3.4}
Let $s\ge r-2$, $r+1, s+2\ge k \ge 1$.
Let $F$ be a $G$-gauge-natural bundle of order $k$.
All natural differential operators
$$
f:C^\infty (\Cla\f M\ucar{\f M} \Lin \f E) \to  C^\infty(F\f E)
$$
which are of order $s$ with respect to classical connections and
of order $r$ with respect to linear connections are of the
form
$$
f(j^s\Lam, j^r K) = g(j^{k-2}\Lam, j^{k-1} K,
\nab^{(k-2,s-1)} R[\Lam], \nab^{(k-1,r-1)} R[K])
$$
where $g$ is a unique natural operator
$$
g:J^{k-2}\Cla \f M\ucar{\f M} J^{k-1} \Lin \f E\ucar{\f M}
C^{(k-2,s-1)}_C\f M\ucar{\f M} C^{(k-1,r-1)}_L\f E
        \to   F\f E\,.
$$
\eTh

\bRm\label{Rm3.5}
From the proof of Theorem \ref{Th3.2} it follows that the operator $g$
is the restriction of a zero order operator defined on the $k$-th
order $G$-gauge-natural bundle
$J^{k-2}\Cla \f M\ucar{\f M} J^{k-1} \Lin \f E\ucar{\f M}
{\M W}^{(k-2,s-1)}\f M\ucar{\f M} {\M U}^{(k-1,r-1)}\f E$.
\eRm

\section{The second $k$-th order reduction theorem for linear
         and classical connections}
\setcounter{equation}{0}
\setcounter{theorem}{0}

Write $(\f E^{p_1,p_2}_{q_1,q_2})_i
\byd\f E^{p_1,p_2}_{q_1,q_2}\ten\ten^iT^*\f M$,
$i\ge 0$, and
set
\beq
(\f E^{p_1,p_2}_{q_1,q_2})^{(k,r)} \byd
        (\f E^{p_1,p_2}_{q_1,q_2})_k\ucar{\f M}\dots\ucar{\f M}
        (\f E^{p_1,p_2}_{q_1,q_2})_r\,,\qquad
        (\f E^{p_1,p_2}_{q_1,q_2})^{(r)} \byd
        (\f E^{p_1,p_2}_{q_1,q_2})^{(0,r)}\,.
\eeq
The $i$-th order covariant differential of sections of
$\f E^{p_1,p_2}_{q_1,q_2}$ with respect to
$(\Lam,K)$ is a natural operator
\beq
\nab^i: C^\infty(\Cla\f M\ucar{\f M}\Lin\f E \ucar{\f M}
        \f E^{p_1,p_2}_{q_1,q_2})
        \to C^\infty((\f E^{p_1,p_2}_{q_1,q_2})_i)
\eeq
which is of order $(i-1)$ with respect to classical and linear
connections and of order $i$ with respect to sections of
$\f E^{p_1,p_2}_{q_1,q_2}$.
Let us note that $\f E^{p_1,p_2}_{q_1,q_2}$
is a $(1,0)$-order $G$-gauge-natural bundle and let us denote by
$V\byd\ten^{p_1} \Rn^n\ten\ten^{q_1}\Rn^{n*}\ten^{p_2}
        \Rn^m\ten\ten^{q_2}\Rn^{m*}$
its standard fiber with coordinates
$(v^A)=(v^{i_1\dots i_{p_1}\lam_1\dots\lam_{p_2}}_{j_1\dots
        j_{q_1}\mu_1\dots\mu_{q_2}})$. By $V_i$ or
$V^{(k,r)}\byd V_k\car \dots \car V_r$, $V^{(r)}\byd V^{(0,r)}$,
we denote
the standard fibers of
$(\f E^{p_1,p_2}_{q_1,q_2})_i$ or
$(\f E^{p_1,p_2}_{q_1,q_2})^{(k,r)}$, respectively.

Hence we have the associated
$W^{(i+1,i+1)}_m G$-equivariant map,
denoted by the same symbol,
\beq
\nab^i: T^{i-1}_m {Q}\car T^{i-1}_m {R} \car
        T^i_m V \to V_i\,.
\eeq

If $(v^A,v^A{}_\lam,\dots,v^A{}_{\lam_1\dots\lam_i})$ are the
induced jet coordinates on $T^i_m V$ (symmetric in all
subscripts) and
$(V^A{}_{\lam_1\dots\lam_i})$ are the canonical coordinates
on $V_i$, then $\nab^i$ is of the form
\begin{align}\label{Eq4.1}
(& V^A{}_{\lam_1\dots\lam_i})\com\nab^i
\\
& =v^A{}_{\lam_1\dots\lam_i} +\pol(T^{i-1}_m {Q}\car
        T^{i-1}_m {R} \car
        T^{i-1}_m V)\,,\nonumber
\end{align}
where $\pol$ is a polynomial on
$T^{i-1}_m {Q}\car T^{i-1}_m {R} \car  T^{i-1}_m V$.

We define the $k$-th order {\em formal Ricci equations},
$k\ge 2$, as follows.
For $k=2$ we have by Remark \ref{Rm2.8}
$$
V^A{}_{[\lam\mu]}-\pol(C^{(0)}_{C}\car C^{(0)}_{{L}}
        \car V) = 0\,.\leqno{(E_2)}
$$
For $k>2$, $(E_k)$ is obtained by the formal covariant
differentiating of $(E_2) - (E_{k-1})$ and antisymmetrization
of the last two formal covariant
differentials. They are of the form
$$
V^A{}_{\lam_1\dots[\lam_i\lam_{i+1}]\dots\lam_k}
        -\pol(C^{(k-2)}_{C}\car C^{(k-2)}_{{L}}
        \car V^{(k-2)}) = 0\,,
\leqno{(E_k)}
$$
$i=1,\dots,k-1$.

\bDf\label{Df4.1}
The $k$-th order {\it formal Ricci subspace} $Z^{(k)}\subset
C^{(k-2)}_{C}\car C^{(k-2)}_{{L}}\car V^{(k)}$ is defined by
equations $(E_2), \dots , (E_k)$, $k\ge 2$. For $k=0,1$
we set $ Z^{(0)}=V$ and
$ Z^{(1)}=V^{(1)}$.
\eDf

In \cite{Jan04} it was proved that $Z^{(k)} $
is a submanifold of $C^{(k-2)}_{C}\car C^{(k-2)}_{{L}}
\car V^{(k)}$
and the restricted morphism
\beq
(\M R^{(k-2)}_{C},\M R^{(k-2)}_{{L}},\nabla^{(k)}):
        T^{k-1}_m {Q}\car T^{k-1}_m {R}
        \car T^k_m V \to Z^{(k)}
\eeq
is a surjective submersion.
Let us consider the projection
$\pr^r_k: Z^{(r)}\to Z^{(k)}$. We have an affine structure on
fibres of the projection $\pr^r_{r-1}: Z^{(r)}\to Z^{(r-1)}$.
It follows from the fact that $Z^{(r)}$ is a subbundle
in $Z^{(r-1)}\car (C_{C,r-2}\car C_{L,r-2} \car V_r)$
given as the space of solutions of
the  system of nonhomogeneous equations ($E_r$).
Let us denote by $Z^{(k,r)}_{z^{(k-1)}}$ the fiber in
$z^{(k-1)}\in Z^{(k-1)}$ of the projection
$\pr^r_{k-1}:Z^{(r)}\to Z^{(k-1)}$.
Then we can consider the fiber product over $Z^{(k-1)}$
$$
 (T^{k-2}_m {Q}\car T^{k-2}_m {R}
        \car T^{k-1}_m V)\ucar{Z^{(k-1)}} Z^{(r)}
$$
and denote it by
$$
T^{k-2}_m {Q}\car T^{k-2}_m {R}
        \car T^{k-1}_m V \car Z^{(k,r)} \,.
$$

\bLm\label{Lm4.2}
If $r+1\ge k\ge 1$, then
the restricted morphism
\begin{multline}
(\pi^{r-1}_{k-2}\car \pi^{r-1}_{k-2}\car \pi^r_{k-1})\car
(\M R^{(k-2,r-2)}_{C},\M R^{(k-2,r-2)}_{{L}},\nabla^{(k,r)}):
\nonumber
\\
      :  T^{r-1}_m {Q}\car T^{r-1}_m {R}
        \car T^r_m V \to T^{k-2}_m Q\car T^{k-2}_m R\car T^{k-1}_m V
\car Z^{(k,r)}\nonumber
\end{multline}
is a surjective submersion.
\eLm

\bPf
The proof of Lemma \ref{Lm4.2} follows from the commutative diagram
\beq
\begin{CD}
        T^{r-1}_m Q\car T^{r-1}_m R\car T^r_m V
        @>(\M R^{(r-2)}_C, \M R^{(r-2)}_L, \nab^{(r)})>>
        Z^{(r)}
\\
        @V{\pi^{r-1}_{k-2}\car \pi^{r-1}_{k-2}\car\pi^r_{k-1}}VV
        @VV{\pr^r_{k-1}}V
\\
        T^{k-2}_m Q \car T^{k-2}_m R \car T^{k-1}_m V
        @>(\M R^{(k-3)}_C,\M R^{(k-3)}_L,\nab^{(k-1)})>>
        Z^{(k-1)}
\end{CD}
\eeq
where all morphisms are surjective submersions.
Hence
\bEq\label{Eq4.2}
(\pi^{r-1}_{k-2}\car\pi^{r-1}_{k-2}\car\pi^r_{k-1})
\car(\M R^{(k-2,r-2)}_C,\M R^{(k-2,r-2)}_L,\nab^{(k,r)})
\eEq
is surjective. For $k=r$ the map $(\M R^{(r-2,r-2)}_C=\M R_{C,r-2},
\M R^{(r-2,r-2)}_L=\M R_{L,r-2},
\nab^{(r,r)}=\nab^r)$ is affine morphism over
$(\M R^{(r-3)}_C,\M R^{(r-3)}_L, \nab^{(r-1)})$
with constant rank, i.e.
$
(\pi^{r-1}_{r-2}\car\pi^{r-1}_{r-2}\car\pi^r_{r-1})
\car(\M R_{C,r-2},\M R_{L,r-2},\nab^{r})
$
is a submersion.
The mapping (\ref{Eq4.2}) is then a composition of surjective submersions.
\ePf

\bTh\label{Th4.3}
Let $S_F$ be a left $W^{(k,k)}_m G$-manifold.
For every $W^{(r+1,r+1)}_m G$-equivariant map
$f: T^{r-1}_m Q\car T^{r-1}_m R\car T^r_m V
\to S_F$ there exists a unique $W^{(k,k)}_m G$-equivariant map
$g: T^{k-2}_m Q\car T^{k-2}_m R\car  T^{k-1}_m V \car Z^{(k,r)}\to S_F$
such that
$$
f = g\com
(\pi^{r-1}_{k-2}\car\pi^{r-1}_{k-2}\car\pi^r_{k-1},
        \M R^{(k-2,r-2)}_C,\M R^{(k-2,r-2)}_L,\nab^{(k,r)})\,.
$$
\eTh

\bPf
Consider the map
\begin{multline}
(\id_{T^{r-1}_m Q}\car\id_{T^{r-1}_m R}\car \pi^{r}_{k-1},\nab^{(k,r)})
        : T^{r-1}_m {Q}\car T^{r-1}_m {R}\car T^{r}_m {V}
        \to     \nonumber
\\
        \to  T^{r-1}_m {Q}\car T^{r-1}_m {R}\car T^{k-1}_m {V}
        \car  V^{(k,r)}  \nonumber
\end{multline}
and denote by $\widetilde{V}^{(k,r)}\subset
        T^{r-1}_m {Q}\car T^{r-1}_m {R}\car T^{k-1}_m {V}
        \car V^{(k,r)}$
its image. By (\ref{Eq4.1}), the restricted morphism
\beq
\widetilde{\nab}^{(k,r)}:T^{r-1}_m {Q}\car
         T^{r-1}_m {R}\car T^{r}_m {V}
        \to
        \widetilde{V}^{(k,r)}
\eeq
is bijective for every $(j^{r-1}_0 \lam,j^{r-1}_0 \gam) \in
T^{r-1}_m {Q}\car  T^{r-1}_m {R}$, so that
$\widetilde{\nab}^{(k,r)}$ is an equivariant diffeomorphism.
Define
\begin{align*}
(\widetilde{\M R}^{(k-2,r-2)}_C,\widetilde{\M R}^{(k-2,r-2)}_L) & :
\widetilde{V}^{(k,r)}
 \to T^{k-2}_m Q\car T^{k-2}_m R\car T^{k-1}_m V\car Z^{(k,r)}
\end{align*}
by
\begin{multline*}
(\widetilde{\M R}^{(k-2,r-2)}_C, \widetilde{\M R}^{(k-2,r-2)}_L)
(j^{r-1}_0 \lam,j^{r-1}_0 \gam,j^{k-1}_0 {\mu},v) =
\\
= (j^{k-2}_0\lam,j^{k-2}_0\gam,
j^{k-1}_0 \mu,{\M R}^{(k-2,r-2)}_C(j^{r-1}_0 \lam),
{\M R}^{(k-2,r-2)}_L(j^{r-1}_0 \lam,j^{r-1}_0 \gam),v)\,,
\end{multline*}
$(j^{r-1}_0 \lam ,j^{r-1}_0 \gam,j^{k-1}_0 {\mu},v)\in
\widetilde{V}^{(k,r)}$.
By Lemma \ref{Lm3.1}
$(\widetilde{\M R}^{(k-2,r-2)}_C,\widetilde{\M R}^{(k-2,r-2)}_L)$
is a surjective submersion.

Thus, Lemma \ref{Lm3.1} and Lemma \ref{Lm3.3} imply that
$(\widetilde{\M R}^{(k-2,r-2)}_C,\widetilde{\M R}^{(k-2,r-2)}_L)$
satisfies
the orbit conditions for the group epimorphism
$\pi^{r+1,r+1}_{k,k}: W^{(r+1,r+1)}_m G\to W^{(k,k)}_m G$ and
there exists a unique
$W^{(k,k)}_m G$-equivariant map
$g: T^{k-2}_m Q\car T^{k-2}_m R\car T^{k-1}_m V \car  Z^{(k,r)}
\to S_F$
such that the diagram
\begin{equation*}\begin{CD}
\widetilde{V}^{(k,r)}
        @>{(\widetilde{\nab}^{(k,r)})^{-1}}>>
        T^{r-1}_m {Q}\car T^{r-1}_m {R}\car T^{r}_m {V}
        @>f>>
        S_F
\\
        @V{(\widetilde{\M R}^{(k-2,r-2)}_C,
        \widetilde{\M R}^{(k-2,r-2)}_L)}VV
        @V{(\pi^{r-1}_{k-2}\car
        \pi^{r-1}_{k-2}\car \pi^{r}_{k-1},
        \M R^{(k-2,r-2)}_C,\M R^{(k-2,r-2)}_L,\nab^{(k,r)})}VV
        @V\id_{S_F}VV
\\
        T^{k-2}_m Q\car T^{k-2}_m R\car T^{k-1}_m V \car  Z^{(k,r)}
        @>{\id}>>
        T^{k-2}_m Q\car T^{k-2}_m R\car T^{k-1}_m V \car Z^{(k,r)}
        @>g>>
        S_F
\end{CD}\end{equation*}
commutes. Hence
$f\com (\widetilde{\nab}^{(k,r)})^{-1}=g\com
(\widetilde{\M R}^{(k-2,r-2)}_C,\widetilde{\M R}^{(k-2,r-2)}_L)$.
Composing both sides with $\widetilde{\nab}^{(k,r)}$,
by considering
$$
(\widetilde{\M R}^{(k-2,r-2)}_C,
\widetilde{\M R}^{(k-2,r-2)}_L)
\com \widetilde{\nab}^{(k,r)}=(\pi^{r-1}_{k-2}\car
\pi^{r-1}_{k-2}\car \pi^{r}_{k-1},
        \M R^{(k-2,r-2)}_C,\M R^{(k-2,r-2)}_L,
\nab^{(k,r)}),
$$
we get
$$
f=g\com (\pi^{r-1}_{k-2}\car \pi^{r-1}_{k-2}\car
\pi^{r}_{k-1},\M R^{(k-2,r-2)}_C,\M R^{(k-2,r-2)}_L,\nab^{(k,r)})
\,.
$$
\vglue-1.3\baselineskip
\ePf

$
T^{k-2}_m Q\car T^{k-2}_m R\car T^{k-1}_m V \car Z^{(k,r)}
$
is closed with respect to the action of the group $W^{(k,k)}_m G$.
The corresponding natural bundle is
$
J^{k-2}\Cla\f M \ucar{\f M} J^{k-2}\Lin\f E \ucar{\f M}
J^{k-1}\f E^{p_1,p_2}_{q_1,q_2}
        \ucar{\f M}Z^{(k,r)}\f E\,.
$
Then the second $k$-order reduction theorem can be
formulated as follows.

\bTh\label{Th4.4}
Let $F$ be a $G$-gauge-natural bundle
of order $k\ge 1$ and let $r+1\ge k$. All
natural differential operators
$f:C^\infty(\Cla\f M\ucar{\f M} \Lin\f E\ucar{\f M}
        \f E^{p_1,p_2}_{q_1,q_2})\to C^\infty(F\f E)$
of order $r$ with respect sections of $\f E^{p_1,p_2}_{q_1,q_2}$
are of the form
\beq
f(j^{r-1}\Lam,j^{r-1}K,j^r\Phi) = g(j^{k-2}\Lam,
        j^{k-2}K, j^{k-1}\Phi,\nab^{(k-2,r-2)}
R[\Lam],\nab^{(k-2,r-2)}
R[K],
        \nab^{(k,r)} \Phi)\,
\eeq
where $g$ is a unique natural operator
$$
g:J^{k-2}\Cla\f M\ucar{\f M} J^{k-2}\Lin\f E \ucar{\f M}
J^{k-1}\f E^{p_1,p_2}_{q_1,q_2}\ucar{\f M}Z^{(k,r)}\f
E\to F\f E\,.
$$
\eTh

\bRm\label{Rm4.5}
The order $(r-1)$ of the above operators with respect to linear and
 classical
connections is the minimal order we have to use. The second reduction
theorem can be easily generalized for any operators of orders $s_1$
or $s_2$ with respect to connections $\Lam$ or $K$,
respectively, where
$s_1\ge s_2-2$, $s_1,s_2\ge r-1$. Then
$$
f(j^{s_1}\Lam,j^{s_2}K,j^r\Phi) =
        g(j^{k-2}\Lam, j^{k-2}K, j^{k-1}\Phi,\nab^{(k-2,s_1-1)}
R[\Lam], \nab^{(k-2,s_2-1)}
R[K],
        \nab^{(k,r)} \Phi)\,.
$$
\eRm

\bRm\label{Rm4.6}
It is easy to see that the second reduction theorem can be generalized
for any number of  fields
$\overset{i}{\Phi}$, $i=1, \dots , m$, of order $(1,0)$
and that any finite order operator
$$
f(j^{s_1}\Lam, j^{s_2}K, j^{r_i}\overset{i}{\Phi}
)\,,\quad s_1,s_2\ge
\max(r_i)-1\,,s_1\ge s_2-2\,,
$$
factorizes through $j^{k-2}\Lam, j^{k-2}K, j^{k-1}\overset{i}{\Phi}$
and sufficiently high covariant differentials of
$R[\Lam]$, $R[K]$, $\overset{i}{\Phi}$.
\eRm



\end{document}